\documentclass[a4paper]{article}
\usepackage{amsfonts}
\usepackage{graphics}

\begin{document}

\bibliographystyle{plain}
\sloppy

\newcommand{\saut}{\par\addvspace{\baselineskip}} 
\newcommand{\demisaut}{\par\vspace{1ex}} 
\newcommand{\eps}{\varepsilon}
\newcommand{\vfi}{\varphi}
\newcommand{\disp}{\displaystyle}
\newcommand{\Int}[1]{\mbox{\rm Int}\left(#1\right)}
\newcommand{\Card}[1]{\mbox{\rm Card}\left(#1\right)}
\newcommand{\diam}[1]{\mbox{\rm diam}\left(#1\right)}
\newcommand{\subsetnoteq}{\subseteq_{\mbox{\rm \hspace{-0.5em}{\scriptsize /}
  \hspace{0.1em}}}}
\newcommand{\qed}{\fbox{\rule{0ex}{0.3ex}\rule{0.3ex}{0ex}}}

\def\IN{\mathbb{N}}
\def\IZ{\mathbb{Z}}

\renewcommand{\labelenumi}{{\rm (\roman{enumi})}}

\newcounter{compteur}[section]
\def\thecompteur{\arabic{section}.\arabic{compteur}}

\newenvironment{resultat}[2][]{\refstepcounter{compteur}\saut%
{\bf \noindent #2\ \thecompteur\ #1 } \it}{\par \saut \rm}

\newenvironment{autre}[2][]{\refstepcounter{compteur}\saut%
{\bf \noindent #2\ \thecompteur\ #1 }}{\par \saut}

\newenvironment{prop}[1][]{\begin{resultat}[#1]{Proposition}}{\end{resultat}}
\newenvironment{lem}[1][]{\begin{resultat}[#1]{Lemma}}{\end{resultat}}
\newenvironment{theo}[1][]{\begin{resultat}[#1]{Theorem}}{\end{resultat}}
\newenvironment{defi}[1][]{\begin{autre}[#1]{Definition}}{\end{autre}}
\newenvironment{rem}[1][]{\begin{autre}[#1]{Remark}}{\end{autre}}
\newenvironment{demo}{\saut\noindent Proof: \par}{\hfill\qed \par\saut}


\title{Mixing $C^r$ maps of the interval without maximal measure
\footnotetext{Israel Journal of Math., {\bf 127}, 253-277, 2002.}}
\author{Sylvie {\sc Ruette}}
\date{Institut de Math\'ematiques de Luminy}

\maketitle

\begin{abstract}
We construct a $C^r$ transformation of the interval (or the torus)
which is topologically  mixing but has no invariant measure of maximal
entropy. 
Whereas the assumption of $C^{\infty}$ ensures existence of
maximal  measures for an interval map, it shows we cannot weaken the
smoothness assumption. 
We also compute the local entropy of the example.
\end{abstract}


\section*{Introduction}
We are interested in topological  dynamical systems on the interval,
that is systems of the form $f\colon I\to I$ where $f$ is
at least continuous and $I$ is a compact interval. 
One can wonder
whether such a system has maximal  measures, i.e. invariant measures
of maximal entropy.

Hofbauer \cite{Hof}, \cite{Hof2} studied piecewise monotone maps, i.e.
interval maps with a finite number of monotone  continuous pieces (the
whole map is not necessarily continuous). 
He proved in this case that
the system admits a non zero finite number of maximal measures if its
topological entropy is positive, and  transitivity implies intrinsic
ergodicity, that is existence of a unique maximal measure. 
For this
purpose, he built a Markov chain which is isomorphic modulo ``small
sets'' with the first system.  
Buzzi \cite{Buz} generalized the
construction of the Markov extension to any continuous interval map.
He showed that the same conclusions as in the piecewise monotone case
hold for $C^{\infty}$ maps. 

One can wonder if these results are still valid under a weaker
regularity assumption, at least in the mixing case. 
Actually, if a
topological dynamical system is expansive and satisfies the
specification property, then it has a unique maximal measure (Bowen
\cite{Bow2}, \cite{Bow3}). 
Specification is a strong property on
periodic points, which must closely follow arbitrary pieces of orbits
(see e.g. \cite{DGS} for more details). 
In the particular case of
continuous interval maps, the system is never expansive, but the
mixing property  implies the specification property (this result is
due to Blokh \cite{Blo1}, see  \cite{Buz3} for the proof). 
More recently, Ruelle \cite{Rue2} worked on positively expansive maps
satisfying specification.

In fact, transitivity is not much weaker than mixing since for any
transitive continuous interval map $f \colon I\to I$ either the map is
mixing or there exist two subintervals $J,K$ such that $J\cup K=I$,
$J\cap K$ is reduced to a single point, $f(J)=K$, $f(K)=J$ and
$f^2|_J, f^2|_K$ are mixing  \cite[p59]{BCop}. 
We also recall that the
topological entropy of any  transitive continuous interval map is
positive (it is greater than or equal to $\frac{\log 2}{2}$
\cite{Blo3}, see  \cite{ALM} for the proof) and, if in addition the
map is Lipschitz, it is finite (this classical result appears in the
proof of Proposition \ref{prop:hG}).

Gurevich and Zargaryan \cite{GZ} built a continuous interval map with
finite entropy which is transitive (in fact mixing) and has no maximal
measure. 
This map has countably many intervals of monotonicity. 
The authors asked is this example can be made smooth on the whole
interval. 
Actually it cannot: the end points $0$ and $1$ are  fixed
points  and the map is not monotone in a neighbourhood of $0$ and $1$;
on the other hand it is not hard to see that a $C^1$ transitive
interval  map must have non zero derivatives at fixed points, hence it
is monotone near these points.

In \cite[Appendix A]{Buz} Buzzi built a $C^r$ interval map which has
no transitive component of maximal entropy, hence it has no maximal
measure. 
He also sketched without details the construction of a $C^r$
interval map with positive entropy which  admits no maximal measure
and which is transitive after restriction to its unique transitive
component (which may be a Cantor set). 
His proof of non existence of
any maximal measure relies on a result of Salama \cite{Sal} 
whose proof turned out to be false (see Theorem 2.3 and Errata in
\cite{Sal2}). 
Nevertheless Buzzi's proof can be modified -- using
extension graphs instead of subgraphs, as we do in Subsection
\ref{subsec:non-existence}~-- so as to be based on another  theorem of
Salama.

The aim of this article is to build for any integer $r\geq 1$ a
$C^r$ mixing interval map which has no maximal measure.  
Transitivity
instead of mixing would be enough, yet it is not more difficult to
prove directly the mixing property. 
This family of examples is
inspired by Buzzi's \cite{Buz}, the important addition is that the system
is transitive on the whole interval.  
Non existence of maximal
measure prevents the metric entropy from being an upper
semi-continuous map on the set of invariant measures.  
This is to be
put in parallel with the result of  Misiurewicz and Szlenk \cite{MS2},
which shows that the topological entropy, considered as a map on the
set of $C^r$ interval transformations, is not upper semi-continuous
for the $C^r$ topology.

In Section 1, we define for any $r\geq 1$ a $C^r$ transformation of
the interval $[0,4]$ which is topologically mixing. 
In fact it is
$C^{\infty}$ everywhere except at one point. 
The map $f_r$ is made of
a countable number of monotone pieces and is Markov with respect to a
countable partition.  
Moreover, it can also be seen as a $C^r$
transformation of the torus by identifying the two end points.  
In the
next section, we study the Markov chain associated with $f_r$  and we
conclude it has no maximal measure, thanks to results of Gurevi\v{c}
\cite{Gur1}, \cite{Gur2}  and Salama \cite{Sal2}. 
As there is an
isomorphism modulo countable sets between the two systems, the
interval map  has no maximal measure either.  
In Section 3, we compute
the local entropy of our examples. 
Buzzi  \cite{Buz} showed that this
quantity bounds the defect in upper-semicontinuity and he gave an
estimate of it depending on the differential order and the spectral 
radius of the derivative. 
Our examples show these bounds are sharp since the two
are realized. 
Moreover, it also equals the topological entropy.  
It
may be of some importance: we conjecture that the Markov extension
admits a maximal measure when the topological entropy is strictly
greater than the local entropy.

In addition to the problem of existence of maximal measure, one can
ask the question of uniqueness of such a measure. 
Recently, Buzzi
\cite{Buz4} proved that, if the interval transformation is
$C^{1+\alpha}$ (i.e. the map is $C^1$ and its derivative is
$\alpha$-H\"older), then there is no measure of positive entropy on
the non Markov part of the system. 
Since a transitive Markov chain
admits at most one maximal measure, a transitive $C^{1+\alpha}$
transformation has a unique maximal measure if it exists.  
For transitive non smooth interval maps we still do not know if
several maximal measures can exist. 
It would imply that the topological
entropy of the critical points would be equal to the topological
entropy of the whole map.

I am indebted to J\'er\^ome Buzzi for many discussions which have led
to the ideas of this paper.

\section{Construction and proof of mixing property}
In this section, we construct a family of $C^r$ maps  $f_r\colon
I\to I$ for $r\geq 1$, where $I=[0,4]$. 
We first give a
general idea of their aspect (see Figure \ref{fig:f}). 
Then we give
some lemmas which will be useful to prove the mixing
property. 
Finally, we define $f_r$ by pieces and check some properties
at each step.  
At the end of the section, the maps $f_r$ are totally
defined and are proved to be mixing.

\subsection{General description}
Let $\lambda\geq 14$ ($\log\lambda$ will be the entropy of $f_r$).
The map $f_r$ is increasing on $[0,1/2]$ and decreasing on $[1/2,1]$. 
Moreover, $f_r(x)=\lambda^r x$ for $0 \leq x
\leq \frac{5}{2} \lambda^{-r},\; f_r(0)=f_r(1)=0,\; f_r(1/2)=4$.

\begin{figure}[htbp]
\begin{center}
\includegraphics{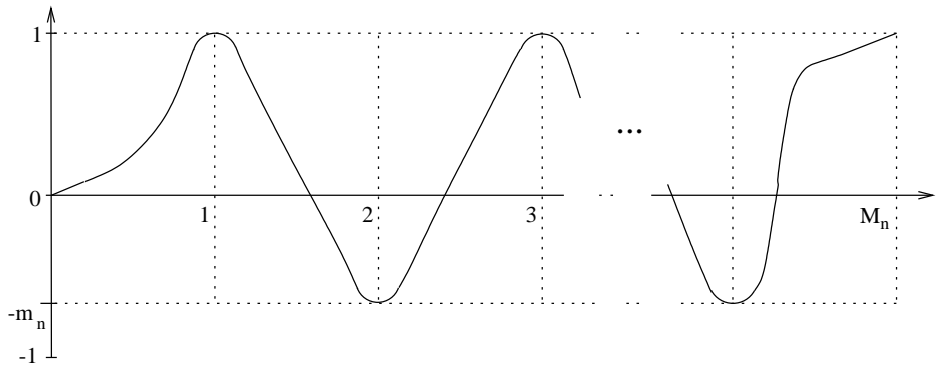}
\caption{the map $s_n$ \label{fig:sn}}
\end{center}
\end{figure}

Let $x_n=1+\frac{1}{n}$ and $y_n=x_n+\frac{1}{2n^2}$ for every $n\geq
1$, and let $M_n$  be a sequence of odd numbers with $(\log
M_n)/n\longrightarrow \log\lambda$.  
We choose a family of
$C^{\infty}$ maps $s_n \colon  [0,M_n] \to [-1,1]$ such
that $s_n$ is nearly 2-periodic and has $M_n$ oscillations; $s_n(0)=0$ and
$s_n(M_n)=1$ (see Figure \ref{fig:sn}).

Then we define $f_r$ on $[x_n,y_n]$ by
$$
f_r(x)=\lambda^{-nr}\left[x_n+(y_n-x_n)s_n\left(M_n\frac{x-x_n}{y_n-x_n}
\right)\right].
$$

In this way, $f_r(x_n)=\lambda^{-nr}x_n$, $f_r(y_n)=\lambda^{-nr}y_n$
and $f_r$ oscillates $M_n$ times between $x_n$ and $y_n$ like $s_n$.
It is worth mentioning that $x_n$ and $y_n$ are periodic points with
period $n+1$, because $f_r$ is linear of slope $\lambda^r$ on
$[0,y_1\lambda^{-r}]$.

On $[y_{n+1},x_n]$, $f_r$ is increasing.

\begin{figure}[bhtp]
\begin{center}
\includegraphics{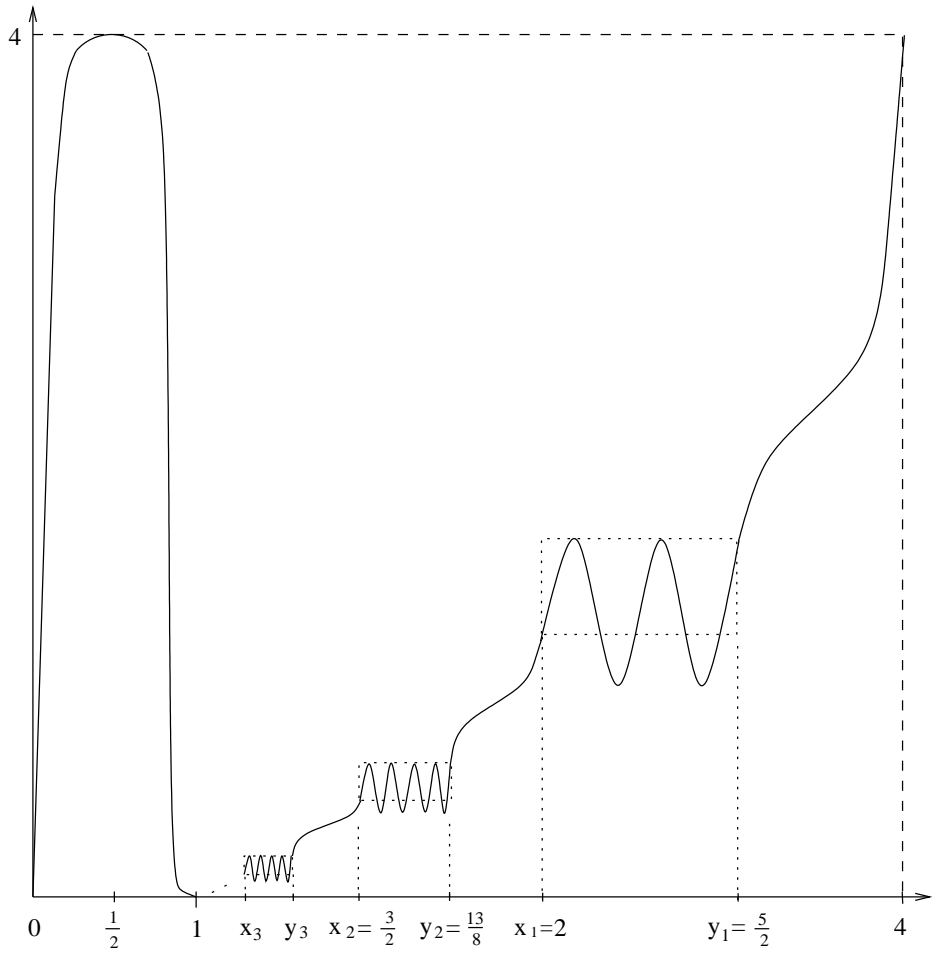}
\caption{the map $f_r$ (scale is not respected) \label{fig:f}}
\end{center}
\end{figure}

Finally, $f_r$ is increasing on $[y_1,4]$, with $f_r(4)=4$. 
Figure \ref{fig:f} gives a general idea of  $f_r$.

\saut The map $f_r$ will be built to be mixing and $C^r$ on $[0,4]$,
and $\|f_r'\|_{\infty}=\lambda^r$. 
Furthermore, the minimum of $s_n$
will be chosen such that $f_r(x)=\lambda^{-nr}y_{n+1}$ if $x$ is a
local minimum of $f_r$ in $]x_n,y_n[$ in order to obtain a Markov map.

This brief description is sufficient to build the Markov chain
associated with $f_r$ and prove that $f_r$ has no maximal measure,
which is done in Section \ref{sec:Markov}. 
The rest of this section, which may be skipped a first reading, is 
devoted to prove that maps satisfying these properties do exist.

\subsection{Method for the proof of mixing property}
We recall the definition of   mixing for a topological dynamical
system.

\begin{defi}
Let $T\colon X\to X$ be a continuous map where $X$ is a compact
metric space.  
The system is (topologically) {\bf mixing} if for every
non empty  open sets $U$ and $V$, there exists $N\geq 0$ such that for
every $n\geq N,\:T^{-n}U\cap V\not= \emptyset$.
\end{defi}

In our case, we will show that for any non degenerate subinterval
$J\subset I$, there exists $n\geq 0$ such that $f_r^n(J)=I$. 
So $f_r^k(J)=I$ for every $k\geq n$ and the system is mixing. 
For this,
we will show that, for some constant $\mu_0>1$, any non degenerate
subinterval $J$ satisfies one of the two following conditions:

\begin{description}
\item[\rm(1)]$\exists k\geq 0$ such that $|f_r^k(J)| \geq \mu_0|J| $,
where $|J|$ denotes the length of $J$,
\end{description}
or
\begin{description}
\item[\rm(2)]$\exists k \geq 0, \exists n\geq 1$ such that either
$0\in f_r^k(J)$ or $\Int{f_r^k(J)}$ contains $x_n$ or $y_n$.
\end{description}

Then it will be enough to show that for any non degenerate
subinterval $J$ containing $0$ or  $x_n$ or $y_n$, there is a $k$ such
that $f_r^k(J)=I$.

Lemma \ref{lem:extremum} says that an  interval near a suitable
extremum satisfies (1) or (2).  
Lemma \ref{lem:point-repulsif}, which
is trivial, says how an interval  containing a repelling periodic
point behaves.

\begin{lem} \label{lem:point-repulsif}
Let $f\colon I\to I$ where $I$ is a compact interval and let $z_0$ be
a periodic point of period $p$. 
Assume $(f^p)'(x)\geq\mu>1$ for every
$x \in[z_0,z_1]$. 
Then for every $x>z_0$ there exists $n\geq 0$ such
that $f^n([z_0,x])\supset [z_0,z_1]$.
\end{lem}

\begin{lem}\label{lem:extremum}
Let $f \colon  I\to I$ be a $C^r$ map where $I$ is a compact interval
and let $z_0$ be an extremum such that $z_1=f^k(z_0)$ is a periodic
point of period $p$. 
Suppose $f^k(x)=z_1+C(x-z_0)^{\alpha}$ for
$|x-z_0|\leq \delta$, with $C\not= 0$ and $\alpha$ an even
integer. 
Let $z_2=f^k(z_0 -\delta)= f^k(z_0+\delta)$.  
Suppose $\disp
f^p|_{[z_1,z_2]}$ is linear of slope $\mu>1$, and
$\frac{\alpha|z_2-z_1|}{\delta}\geq \mu_0$.  
Then for every non
degenerate interval $J\subset [z_0-\delta,z_0+\delta]$, there exists
$n\geq 0$ for which  one of the following cases holds:
\begin{enumerate}
\item $|f^n(J)|\geq \mu_0 |J|$.
\item $z_2\in \Int{f^n(J)}$.
\end{enumerate}
\end{lem}

\begin{demo}
Let $J=[a,b]$ be an interval in $[z_0-\delta,z_0+\delta]$ with $a<b$.
If $z_0\in J$ then $f^k(J)=[z_1,y]$ for some $y$. 
The hypotheses imply
that $f^p(z_2)>z_2$, hence $z_2$ cannot be an end point of $I$ and one
can choose $1<\mu'<\mu$ and  $z_3>z_2$ such that $(f^p)'(x)>\mu'$ for
all $x\in [z_1,z_3]$. 
According to Lemma \ref{lem:point-repulsif}
there exists $n$ such that $f^n(J)\supset  [z_1,z_3]$, thus $z_2\in
\Int{f^n(J)}$, which is (ii).

Now assume that $z_0\not\in J$. 
We restrict to the case $C>0$ and
$z_0<a<b \leq z_0+\delta$.  
Let $J'=f^k(J)=[a',b']\subset]z_1,z_2]$
and $g=f^p$.  
The point $z_1$ is fixed for $g$ and $g$ is linear of
slope $\mu>1$ on $[z_1,z_2]$, so the map $g$ can be iterated  on $J'$
as long as $g^m(b')\leq z_2$. 
Let $m$ be the first integer satisfying
$g^m(b')>z_2$.  
Then there are two cases:
\begin{itemize}
\item $g^m(a')<z_2<g^m(b')$, which implies (ii) with $n=mp+k$.
\item $z_2\leq g^m(a')<g^m(b')$.

In this case, as $(f^k)'$ is positive and increasing on
$[z_0,z_0+\delta]$ one gets $ |J'|\geq \alpha C(a-z_0)^{\alpha-1}|J|$
and
$$
|f^{mp+k}(J)|\geq \mu^m \alpha C (a-z_0)^{\alpha-1}|J|.
$$ 
But $g^m(a')-z_1=\mu^m(a'-z_1)\geq z_2-z_1$, so
$$
\disp \mu^m\geq
\frac{z_2-z_1}{a'-z_1}=\frac{z_2-z_1}{C(a-z_0)^{\alpha}},
$$ 
and
$$
\disp |f^{mp+k}(J)|\geq \frac{\alpha |z_2-z_1| }{ |a-z_0| } |J| \geq
\frac{\alpha |z_2-z_1|}{\delta}|J|\geq \mu_0 |J|.
$$
\end{itemize}
\end{demo}

We add a lemma which will be useful for some estimates.

\begin{lem} \label{lem:lambda-n}
Let $\lambda \geq 8$ and $[\cdot]$ refer to the entire part of a
number.  
Then for all $n\geq 1$:
\begin{enumerate}
\item 
$\disp \frac{\lambda^n}{n^2}\geq \lambda $.
\item
$\disp \frac{\lambda^n}{2n^2}\leq
2\left[\frac{\lambda^n}{2n^2}\right]-1  \leq \frac{\lambda^n}{n^2} $.
\item
$\disp 2\left[\frac{\lambda^n}{2n^2}\right] -1 \geq \lambda -3$.
\end{enumerate}
\end{lem}

\begin{demo}
(i) is obtained by studying  the  function $x\mapsto
\lambda^{x-1}-x^2$.

For the first inequality of (ii), we  write
$$
2\left[\frac{\lambda^n}{2n^2} \right]-1\geq
\frac{\lambda^n}{2n^2}+\left(\frac{\lambda^n}{2n^2}-3\right) \geq
\frac{\lambda^n}{2n^2}
$$ 
thanks to (i).  
The second inequality is obvious.

(iii) comes from $\disp 2\left[\frac{\lambda^n}{2n^2} \right]-1\geq
\frac{\lambda^n}{n^2}-3$ and from (i).
\end{demo}
 
\subsection{Construction of $f_r$ on $[1,y_1]$} \label{subsec:1.3}
Recall that $\lambda\geq 14$, $f_r(1)=0$, $x_n=1+\frac{1}{n}$ and
$y_n=x_n+\frac{1}{2n^2}\,$ for $n\geq 1$; in particular
$y_1=\frac{5}{2}$.  
In this subsection, we define $f_r$ on $[1,y_1]$
with more details.  
For this purpose, we define $f_r$ on each
$[x_n,y_n]$ and then on each  $[y_{n+1},x_n]$. 
At each step, we check
that the various pieces can be  glued together in  a $C^{\infty}$ way
and  $|f_r'(x)|\leq \lambda^r$ for $x\in [1,y_1]$.  
In addition, we
show that $f_r$ is $C^r$ on the right of $1$. 
Finally, we focus on the
mixing property.  
The map $f_r$ is not totally defined yet, but at
this stage we only need to  know that  $f_r(x)=\lambda^r x$ for $0\leq
x\leq \frac{5}{2}\lambda^{-r}$ and  $f_r(\frac{1}{2})=4$ in order to
prove that any non degenerate subinterval of $[1,y_1]$ satisfies (1)
or (2) with $\mu_0=\frac{4}{3}$. 
Then we show that for an open
interval $J$ containing $x_n$ or $y_n$ there is a $k$ satisfying
$f_r^k(J)=[0,4]$.

\subsubsection{On the subintervals $[x_n,y_n]$} \label{subsubsec:1.3.1}
Set  $ \disp M_n=2 \left[\frac{\lambda^n}{2n^2}\right]-1\mbox{ (where
$[\cdot]$  denotes the entire part)}$, $\disp
m_n=1-\frac{1}{(n+1)^2}$,   $\delta=\lambda^{-r}$, $
C=\frac{1}{4\delta^2}$ and  $k_n=\frac{2\lambda^r}{M_n}.$

First, we choose a sequence of $C^{\infty}$ functions $s_n\colon
[0,M_n]\to [-m_n,1]$ satisfying:

\begin{description}
\item[\rm(3)] $s_n(0)=0,\: s_n(M_n)=1$, $s_n$ is increasing on each
$[2k,2k+1] \:(0\leq k \leq (M_n-1)/2)$, $s_n$ is decreasing on each
$[2k+1,2k+2] \:(0\leq k \leq (M_n-3)/2)$.
\item[\rm(4)] $s_n(x)=1-C(x-a)^2$ for $|x-a|\leq \delta$ if $a$ is a
local maximum of $s_n,\, a\not= M_n$, and $s_n(x)=-m_n+C(x-b)^2$ for
$|x-b|\leq  \delta$ if $b$ is a local minimum of $s_n,\; b\not= 0$.
\item[\rm(5)] $s_n(x)=k_n(x-M_n)+1$ for $M_n-\delta\leq x \leq M_n$.
\item[\rm(6)] $s_n(x)=k_n x$ for $x\in [0,\delta]$.
\item[\rm(7)]  $\forall k\geq 1,\:\exists A_k,\:\forall n\geq 1,\:
\|s_n^{(k)}\|_{\infty}\leq A_k$.
\item[\rm(8)] $\|s_n'\|_{\infty}\leq \lambda^r$ and $|s_n'(x)|\geq
\min\{1/2,k_n\}$ if $|x-d|\geq \delta$ for all local  extrema $d\in
]0,M_n[$.
\end{description}

Property (7) can be fulfilled because $m_n$ and $k_n$ are bounded
($3/4\leq m_n\leq 1$,  $k_n\leq \lambda^r$)  and the maps $s_n$ have a
2-periodic looking.

If $d$ is a local extremum in $]0,M_n[$, then $|s_n(d-\delta)-s_n(d)|=
|s_n(d+\delta)-s_n(d)|=1/4$; moreover $|s_n(\delta)-s_n(0)|\leq 1/4$
and $|s_n(M_n-\delta)-s_n(M_n)|\leq 1/4$. 
Thus if $d$ and $d'$ are two
successive extrema in $[0,M_n]$ the absolute value of the average
slope between $d+\delta$ and $d'-\delta$ is at least
$\disp\frac{m_n+1/2}{1-2\delta} >\frac{1}{2}$ and is less that
$2$. 
Since $|s_n'(d+\delta)|=|s_n'(d-\delta)|= \frac{\lambda^r}{2}$
for any extremum $d\in ]0,M_n[$, Property (8) can be fulfilled.

\saut Secondly, recall that $f_r$ is defined for $x\in [x_n,y_n]$ by
$$ 
f_r(x)=\lambda^{-nr}\left[x_n+(y_n-x_n)s_n\left(M_n\frac
{x-x_n}{y_n-x_n}\right)\right].
$$ 
Now, we look at the $C^r$ character
of $f_r$ near $1$. 
The definition  of $f_r$ gives
$$
f_r^{(k)}(x)\,=\,\frac{\lambda^{-nr}M_n^k}{(y_n-x_n)^{k-1}}\,s_n^{(k)}
\left(M_n\frac{x-x_n}{y_n-x_n}\right) \mbox{ for } x\in [x_n,y_n],
$$
where $f_r^{(k)}(x_n)$ and $f_r^{(k)}(y_n)$ are to be understood as
left (resp. right) derivatives at this stage.

Since $M_n\leq \frac{\lambda^n}{n^2}$,  Property (7) leads to
$$
|f_r^{(k)}(x)|\leq \lambda^{-n(r-k)}n^{-2}2^{k-1} A_k.
$$
One has $A_1=\lambda^r$ by (8), thus $\,|f'(x)|\leq \lambda^r$. 
Moreover, for
$0\leq k\leq r$,
$$
|f_r^{(k)}(x)|\to 0 \mbox{ when } x\to 1,\;x\in\bigcup_{n\geq 1}
[x_n,y_n].
$$ 
Notice that the main factor in this estimate is
$\lambda^{-n(r-k)}$. 
If $k>r$, the $k$-th derivative $f_r^{(k)}$ does
not tend to zero any longer and it can be shown that $f_r$ cannot be
$C^{r+1}$ at point $1$.

\saut As $f_r(x)=\lambda^r x$ for $x\in[0,y_1\lambda^{-r}]$, the
$(n+1)$-th iterate of the map on $[x_n,y_n]$ is given by
$f_r^{n+1}(x)=\lambda^{nr}f_r(x)$.
 
Notice that $m_n$ is chosen such that $\min\{f_r^{n+1}(x)\colon
x\in[x_n,y_n]\}=y_{n+1}$. 
Moreover $f_r^{n+1}(x_n) =x_n$ and
$f_r^{n+1}(y_n)=y_n$.
 
\saut We sum up the previous results in two lemmas, the first one is
about derivatives and the second summarizes the behaviour of $f_r$ on
$[x_n,y_n]$.

\begin{lem}\label{lem:Cr1}
\begin{itemize}
\item $|f_r'(x)|\leq \lambda^r$ for $x\in [x_n,y_n]$.
\item
$\disp\lim_{\mbox{\tiny \begin{tabular}{c} {\scriptsize $x\to 1$} \\
$x\in \cup_{n\geq 1}[x_n,y_n]$ \end{tabular}}} f_r^{(k)}(x)=0 \mbox{
for }\, 0\leq k\leq r.$
\end{itemize}
\end{lem}

\begin{lem}\label{lem:xn-yn}
Let $\disp t_i^n=x_n+\frac{i(y_n-x_n)}{M_n}$ for
$i=0,\cdots,M_n$. 
Then
\begin{itemize}
\item $f_r$ is monotone on $[t_{i-1}^n,t_i^n]$, $1\leq i\leq
M_n$.
\item $f_r(t_i^n)=\lambda^{-nr}y_{n+1}$ if $i$ is even, $i\not=0$, and
$f_r(x_n)=\lambda^{-nr}x_n$.
\item $f_r(t_i^n)=\lambda^{-nr}y_n$ if $i$ is odd.
\end{itemize}
\end{lem}

\subsubsection{On the subintervals $[y_{n+1},x_n]$} \label{subsubsec:1.3.2}
We define
$$
w_n=y_{n+1}+\frac{n+2}{2n(n+1)^2 M_{n+1}k_{n+1}}.
$$
We have $w_n\in]y_{n+1},x_n[$. 
On $[y_{n+1},w_n]$, we define $f_r$ to
be  affine of slope $\lambda^{-(n+1)r}M_{n+1}k_{n+1}$ (recall that
$f_r(y_{n+1})=\lambda^{-(n+1)r}y_{n+1}$ is already defined).  
Because
of this definition $f_r$ is affine (thus $C^{\infty}$) in  a
neighbourhood of $y_{n+1}$. 
Moreover
$$
f_r(w_n)= \lambda^{-(n+1)r}y_{n+1}+
\lambda^{-(n+1)r}\frac{n+2}{2n(n+1)^2}=\lambda^{-(n+1)r}\left(1+\frac{1}{n}
\right),
$$  
so $f_r(w_n)=\lambda^{-(n+1)r}x_n$ and
$f_r^{n+2}(w_n)=x_n$.  
As we are going to extend $f_r$ in a
$C^{\infty}$ way on $[w_n, x_n]$, we will have
$$
f_r'(x_n)=2\lambda^{-(n-1)r},\ f_r'(w_n)=2\lambda^{-nr}, \mbox{ and
}f_r^{(k)}(x_n)=f_r^{(k)}(w_n)=0 \mbox{ for }k\geq 2.
$$

Set $h_n=f_r(x_n)-f_r(w_n)$ and $l_n=x_n-w_n$. 
We compute upper and
lower bounds for $h_n$ and $l_n$.  
First
$$
h_n=\lambda^{-nr}\left(x_n-\lambda^{-r}y_{n+1}-\lambda^{-r}
\frac{n+2}{2n(n+1)^2}\right)\leq 2\lambda^{-nr}.
$$ 
We have
$$
\frac{n+2}{2n(n+1)^2}=\frac{3}{8}\ \mbox{ for }\ n=1
$$ 
and
$$
\frac{n+2}{2n(n+1)^2}=\frac{n^2+2n}{2n^2(n+1)^2}\leq
\frac{1}{2n^2}\leq  \frac{1}{8}\ \mbox{ for }\ n\geq 2.
$$

Since $x_n\geq 1, y_{n+1}\leq y_2=\frac{13}{8}$ and
$\frac{n+2}{2n(n+1)^2}\leq \frac{3}{8}$ for all $n\geq 1$ one gets
$$
h_n\geq \lambda^{-nr}(1-2\lambda^{-r})\geq \frac{6}{7}\lambda^{-nr}.
$$

For $l_n$ one has
\begin{eqnarray*}
l_n &=& 1+\frac{1}{n}-1-\frac{1}{n+1}-\frac{1}{2(n+1)^2}-
\frac{n+2}{2n(n+1)^2M_{n+1}k_{n+1}}\\ &=&
\frac{n+2}{2n(n+1)^2}\left(1- \frac{1}{2\lambda^r}\right).
\end{eqnarray*}
As $\frac{n+2}{2n(n+1)^2}\leq \frac{3}{8}$, one gets  $l_n\leq
\frac{3}{8}$ too. 
Moreover
$$
\frac{n+2}{2n(n+1)^2}\geq \frac{1}{2(n+1)^2}
$$  
and
$\frac{1}{2\lambda^r}\leq \frac{1}{2}$ thus $l_n\geq
\frac{1}{4(n+1)^2}$.  
Finally we obtain the inequalities
$$
\frac{6}{7}\lambda^{-nr}\leq h_n\leq2\lambda^{-nr}
\hspace{2em} \mbox{ and } \hspace{2em} \frac{1}{4(n+1)^2}\leq l_n\leq
\frac{3}{8}.
$$

We normalize $f_r$ on $[w_n,x_n]$ as follows: we define
$\vfi_n\colon[0,1]\to [0,1]$ by
$$
\vfi_n(x)=h_n^{-1}[f_r(w_n+l_n x) -f_r(w_n)].
$$ 
The aim of this normalization is to check that the sequence  $\vfi_n$
can be chosen with  uniformly bounded $k$-th derivatives then to come
back to  $f_r$ and show that $f_r$ is $C^r$ at the right of $1$.  
We want to have
$$
\vfi_n'(1)=2h_n^{-1}l_n\lambda^{-(n-1)r},\
\vfi_n'(0)=2h_n^{-1}l_n\lambda^{-nr},
$$  
and
$$
\vfi_n^{(k)}(0)= \vfi_n^{(k)}(1)=0 \mbox{ for }k\geq 2,
$$ 
thus
$\vfi_n'(1)\leq \frac{7}{8}\lambda^r$ and $\vfi_n'(0)\leq
\frac{7}{8}$.  
Consequently, it is possible to build a sequence of
functions $\vfi_n$ satisfying these conditions and the following
additional conditions:
$$
\forall k\geq 1,\,\exists B_k,\, \forall n\geq 1,\,
\|\vfi_n^{(k)}\|_{\infty}\leq B_k
$$  
and
$$
\forall x\in[0,1],\:  \frac{2}{3}\vfi_n'(0)\leq\vfi_n'(x)\leq
\lambda^r.
$$

By definition of $\vfi_n$, the derivatives of $f_r$ are given by
$$ 
f_r^{(k)}(x)=h_n l_n^{-k}\vfi_n^{(k)}\left(\frac{x-w_n}{l_n}\right)
\leq \lambda^{-nr}(n+1)^{2k}2^{2k+1}B_k \mbox{ for } w_n\leq x\leq
x_n,
$$  
hence for every $k\geq 0$
$$
f_r^{(k)}(x)\to 0\mbox{ when } x\to 1, x\in\bigcup_{n\geq
1}[y_{n+1},x_n].
$$  
Moreover, $\frac{4}{3}\lambda^{-nr}\leq
f_r'(x)\leq \lambda^r h_nl_n^{-1}$  for every $x\in [w_n,x_n]$ and
$$
h_nl_n^{-1}\leq\frac{8(n+1)^2}{\lambda^{nr}}\leq 1 \mbox{ by Lemma
\ref{lem:lambda-n} (i)}.
$$

The next lemma recalls the behaviour of $f_r$ on $[y_{n+1},x_n]$.

\begin{lem}\label{lem:yn+1-xn} 
\begin{itemize}
\item $\frac{4}{3}\lambda^{-nr}\leq f_r'(x) \leq \lambda^r$ for $x\in
[y_{n+1},x_n]$.
\item $f_r(w_n)=\lambda^{-(n+1)r}x_n$.
\item
$\disp\lim_{\mbox{\tiny \begin{tabular}{c} {\scriptsize $x\to 1$} \\
$x\in \bigcup_{n\geq 1}[y_{n+1},x_n]$ \end{tabular}}} f_r^{(k)}(x)=0
\mbox{ for }\, 0\leq k\leq r.$
\end{itemize}
\end{lem}

\subsubsection{Beginning of the proof of the mixing property}
We show that any non degenerate subinterval $J\subset[1,y_1]$
satisfies (1) or (2) with $\mu_0=\frac{4}{3}$. 
It is sufficient to
consider $J\subset [x_n,y_n]$ or $J\subset [y_{n+1},x_n]$.

First, we look at $[y_{n+1},x_n]$. 
For $x\in[y_{n+1},x_n]$,
$f_r^{n+1}(x)=\lambda^{nr}f_r(x)$ and the derivative of $f_r$
satisfies $f_r'(x)\geq \frac{4}{3}\lambda^{-nr}$ by Lemma
\ref{lem:yn+1-xn}, so $|f_r^{n+1}(J)|\geq \frac{4}{3} |J|$ if
$J\subset [y_{n+1},x_n]$.

Now, we focus on $[x_n,y_n]$.  
According to Property (8),
$s_n'(x)\geq\min\{k_n,1/2\}$ for all $x\in [M_n-1+\delta,M_n]$ thus
$$
(f_r^{n+1})'(x)\geq\min\{M_nk_n,M_n/2\}\geq 2\mbox{ for all }
x\in\left[y_n-\frac{(y_n-x_n)(1-\delta)}{M_n},y_n\right].
$$ 
Because of Property (4), $s_n(M_n-1+\delta)= -m_n+1/4<0$, thus
$$
f_r^{n+1}\left(y_n-\frac{(y_n-x_n)(1-\delta)}{M_n}\right)<x_n.
$$
Let $\disp t_n=\frac{y_n-x_n}{\lambda^rM_n}$, then according to Lemma
\ref{lem:point-repulsif},  there exists an integer $\alpha$ such that
$f_r^{(n+1)\alpha}([y_n-t_n,y_n])\supset  [x_n,y_n]$, so there exists
$z\in [y_n-t_n,y_n[$ with  $f_r^{(n+1)\alpha}(z)=x_n$.  
Because of the
choice of $t_n$ and Property (5), $f_r^{n+1}$ is affine of slope
$k_nM_n=2\lambda^r$ on $[y_n-t_n,y_n]$. 
Let $k\geq 0$ be the maximal
integer  $i$ such that $\lambda^{ri}(y_n-z)\leq t_n$. 
Then
$z_n=y_n-\lambda^{rk}(y_n-z)$ belongs to
$\left[y_n-t_n,y_n-\frac{t_n}{2\lambda^r}\right]$ and
$f_r^{(n+1)\alpha_n}(z_n)=x_n$ if $\alpha_n=\alpha+k$.

\saut Set $\disp \delta_n=\sqrt{\frac{(y_n-z_n)(y_n-x_n)}{CM_n^2}}$,
and let $a$ be a local maximum of $f_r$ on $]x_n,y_n[$.  
If $|t|\leq
\delta_n$, then
$$
\left|\frac{M_nt} {y_n-x_n}\right|^2\leq\frac{y_n-z_n}{C(y_n-x_n)}\leq
\frac{1}{CM_n\lambda^r} \leq\delta^2.
$$ 
Now we check the hypotheses of Lemma \ref{lem:extremum} for the
extremum $a$:
\begin{itemize}
\item $f_r^{n+1}(a)=y_n$ and $f_r^{n+1}(y_n)=y_n$.
\item $\disp f_r^{n+1}(a+t)=y_n-\frac{C(M_nt)^2}{y_n-x_n}$ if $|t|\leq
\delta_n$ (because of Property (4)).
\item $f_r^{n+1}(a-\delta_n)=f_r^{n+1}(a+\delta_n)=z_n$.
\item $f_r^{n+1}$ is linear on $[z_n,y_n]$, with a slope $k_nM_n\geq
2$.
\item $\disp
\frac{2|z_n-y_n|}{\delta_n}=2\sqrt{\frac{CM_n^2(y_n-z_n)}{y_n-x_n}}\geq
2\sqrt{\frac{CM_n^2
t_n}{2\lambda^r(y_n-x_n)}}=2\sqrt{\frac{M_n}{8}}$\\ and the last
quantity is greater than $2$ because $M_n\geq \lambda-3$ by Lemma
\ref{lem:lambda-n} (iii) and $\lambda\geq 14$.
\end{itemize}
Consequently, we can apply Lemma \ref{lem:extremum} at this maximum:
for any  non degenerate subinterval $J\subset
[a-\delta_n,a+\delta_n]$, there  exists $k$ such that either
$z_n\in\Int{f_r^k(J)}$ or $|f_r^k(J)|\geq 2|J|$. 
Since
$f_r^{(n+1)\alpha_n}(z_n)=x_n$ and $f_r^{(n+1)\alpha_n}$ is a local
homeomorphism in a neighbourhood of $z_n$, if $z_n\in\Int{f_r^k(J)}$
then $x_n\in\Int{f_r^{k'}(J)}$ with $k'=k+(n+1)\alpha_n$.

\saut Set $\disp
\delta_n'=\sqrt{\frac{(w_n-y_{n+1})(y_n-x_n)}{CM_n^2}}$ and let $b$ be
a local minimum of $f_r$ on $]x_n,y_n[$. 
If $|t|\leq  \delta_n'$, then
$$
\left|\frac{M_nt} {y_n-x_n}\right|^2\leq
\frac{w_n-y_{n+1}}{C(y_n-x_n)}=
\frac{2n(n+2)}{(n+1)^2\lambda^{3r}}\leq \frac{2}{\lambda^{3r}}\leq
\delta^2.
$$

We check the hypotheses of Lemma \ref{lem:extremum} for the extremum b:
\begin{itemize}
\item $f_r^{n+1}(b)=y_{n+1}$ and $f_r^{n+2}(y_{n+1})=y_{n+1}$.
\item $\disp f_r^{n+1}(b+t)=y_{n+1}+\frac{C(M_nt)^2}{y_n-x_n}$ if
$|t|\leq  \delta_n'$ (because of Property (4)).
\item $f_r^{n+1}(b-\delta_n')=f_r^{n+1}(b+\delta_n')=w_n$ and
$f_r^{n+2}(w_n)=x_n$.
\item $f_r^{n+2}$ is linear on $[y_{n+1},w_n]$ of slope
$M_{n+1}k_{n+1}\geq 2$.
\item $\disp \frac{2|y_{n+1}-w_n|}{\delta_n'}\geq 2$.
\end{itemize}
To prove the last point, define
$$
C_n = \left(\frac{w_n-y_{n+1}}{\delta_n'}\right)^2
=\frac{n(n+2)M_n^2\lambda^r}{8(n+1)^2}.
$$
One has $M_n\geq \lambda-3$ (Lemma \ref{lem:lambda-n} (iii)),
$\lambda\geq 14$ and
$$
\frac{2n(n+2)}{(n+1)^2}=\frac{(n+1)^2+n^2+2n-1}{(n+1)^2}>1,
$$
thus $C_n\geq \frac{14\times 11^2}{16}>1$.

Hence we can apply Lemma \ref{lem:extremum} to this extremum: for any
non degenerate subinterval $J\subset [b-\delta_n',b+\delta_n']$, there
exists $k$ such that either $x_n\in\Int{f_r^k(J)}$ or $|f_r^k(J)|\geq
2|J|$.

\saut If $\disp |x-d|\geq \delta|y_n-x_n|/M_n$ for all  local extrema
$d\in]x_n,y_n[$,  then  $|(f_r^{n+1})'(x)|\geq
\min\{2\lambda^r,M_n/2\}\geq 2$ according to Property (8).  
If $a\in
]x_n,y_n[$ is a  local maximum and $\disp \delta_n\leq |x-a|\leq
\frac{\delta|y_n-x_n|}{M_n}$, then
$$
|(f_r^{n+1})'(x)|\geq|(f_r^{n+1})'(a+\delta_n)|= \frac{2M_n^2
C\delta_n}{y_n-x_n} \geq \sqrt{M_n/2}.
$$
If $b \in ]x_n,y_n[$ is a  local minimum and  $\disp\delta_n'\leq
|x-b|\leq \frac{\delta|y_n-x_n|}{M_n}$, then
$$
|(f_r^{n+1})'(x)|\geq|(f_r^{n+1})'(b+\delta_n')|= \frac{2M_n^2
C\delta_n'}{y_n-x_n}=\lambda^{\frac{r}{2}}\frac{M_n}{2}
\sqrt{\frac{2n(n+2)}{(n+1)^2}}\geq \lambda^{\frac{r}{2}}M_n/2.
$$ 
Consequently, $|(f_r^{n+1})'(x)|\geq 2$ if  for all local maxima $a$,
$|x-a|\geq\delta_n$ and for all local minima $b$,
$|x-b|\geq\delta'_n$.

Finally,  if $J$ is a  non degenerate subinterval of $[x_n,y_n]$,
there exists $k$ such that either $|f_r^k(J)|\geq 2|J|$ or
$\Int{f_r^k(J)}$ contains $x_n$.  
Together with the previous result on
$[y_{n+1},x_n]$ it gives:

\begin{lem} \label{lem:lem1.1}
If $J$ is a non degenerate subinterval of $[1,y_1]$, there exist
$k\geq 0$ and $n\geq 1$  such that either $|f_r^k(J)|\geq \frac{4}{3}|J|$ or
$x_n\in \Int{f_r^k(J)}$ or $y_n\in \Int{f_r^k(J)}$.
\end{lem}

The point $x_n$ is periodic of period $n+1$, and $(f_r^{n+1})'(x)\geq
2$  for $ x_n\leq x \leq x_n +\frac{y_n-x_n}{2M_n}$.  
In this
situation, we can apply Lemma \ref{lem:point-repulsif}. 
For any
interval $J=[x_n,x]$ with $x>x_n$ there exists $k$ such that
$f_r^k(J)\supset [x_n,x_n+\frac{y_n-x_n}{2M_n}]$.  
But
$$
f_r^{n+1}\left(x_n+\frac{y_n-x_n}{2M_n}\right)\geq
x_n+\frac{y_n-x_n}{M_n} \ \mbox{ and }\
f_r^{n+1}\left(x_n+\frac{y_n-x_n}{M_n}\right)=y_n.
$$ 
Hence $f_r^{k+2(n+1)}(J)\supset [x_n,y_n]$.

We  do the same thing for $y_n$: for any interval $J=[y,y_n]$ with
$y<y_n$ there exists $k$ such that $f_r^k(J)\supset[x_n,y_n]$.

Moreover
$$
f_r^{2(n+1)}([x_n,y_n])=f_r^{n+1}([y_{n+1},y_n])=
[\lambda^{-1}y_{n+1}, y_n]\supset [1/2,1],
$$ 
so $f_r^{2(n+1)+1}([x_n,y_n])=[0,4]$.  
This leads to the next lemma.

\begin{lem} \label{lem:lem1.2}
If $J$ is an open subinterval with $x_n\in J$ or $y_n\in J$, then
there exists $k\geq 0$ such that $f_r^k(J)=[0,4]$.
\end{lem}

\subsection{Construction of $f_r$ on $[0,1]$ and $[y_1,4]$ and end of the 
proof of the mixing property} \label{subsec:1.4}

Recall that $f_r(x)=\lambda^r x$ for $0\leq x\leq \frac{5}{2}
\lambda^{-r}$ and $\delta=\lambda^{-r}$. 
We define $f_r$ near the
points $1/2,\:1$ and $4$ as follows:
\begin{itemize}
\item $f_r(x)=4-C_0(x-1/2)^2$ for $|x-1/2|\leq\delta$, with
$C_0=\frac{3}{2}\delta^{-1}$.
\item $f_r(x)=C_1(x-1)^{\alpha_1}$ for $1-\delta\leq x\leq 1$, with
$\alpha_1=2r$ and $C_1=\delta^{1-\alpha_1}$.
\item $f(x)=4+\lambda^r(x-4)$ for $4-\frac{3}{2}\delta\leq x\leq 4$.
\end{itemize}

The definition of $f_r$ on the left of $1$, together with Lemmas
\ref{lem:Cr1} and \ref{lem:yn+1-xn}, leads to the next lemma.

\begin{lem}
$f_r$ is $C^r$ in a neighbourhood of $1$.
\end{lem}

Now we complete the map such that the  pieces are glued together in a
$C^{\infty}$ way (except at $1$ where $f_r$ is only $C^r$).  
As $f_r'(1/2-\delta)=3$ and
$$
\frac{f_r(1/2-\delta)-f_r(\frac{5}{2}
\delta)}{(1/2-\delta)-\frac{5}{2}\delta}=
\frac{3-3\lambda^{-r}}{1-7\lambda^{-r}} \in[2,6],
$$ 
the map can be chosen such that $3/2\leq f_r'(x)\leq \lambda^r$ for
every $x\in [\frac{5}{2}\delta,\frac{1}{2}-\delta]$.  
In the same way,
it is possible to have $-\lambda^r\leq f_r'(x)\leq -3/2$  for every
$x\in [1/2+\delta,1-\delta]$ because $f_r'(1/2+\delta)=-3$,
$f_r'(1-\delta) =-2r$ and
$$
\frac{f_r(1/2+\delta)-f_r(1-\delta)}{1/2-2\delta}=
\frac{8-5\lambda^{-r}}{1-4\lambda^{-r}}\in [7,12].
$$
Finally, $f_r'(y_1)=2$ because of the earlier  construction of $f_r$
on $[x_1,y_1]$ (see parag. \ref{subsubsec:1.3.1})  and
$$
\frac{f_r(4-\frac{3}{2}\delta)-
f_r(y_1)}{(4-\frac{3}{2}\delta)-y_1}=\frac{4-4\lambda^{-r}}{\frac{3}{2}-\frac{3}{2}\lambda^{-r}}=
\frac{8}{3}.
$$
Hence it is possible to have $\frac{3}{2} \leq f_r'(x)\leq \lambda^r$
for $y_1\leq x\leq 4$.

\saut Consequently, $\frac{3}{2}\leq |f_r'(x)|\leq \lambda^r$ if $x\in
[0,\frac{1}{2} -\delta]\cup[\frac{1}{2}+\delta,1-\delta]\cup[y_1,4]$.

\saut A quick check shows that Lemma \ref{lem:extremum} can be applied
to the two extrema $1/2$ and $1$ (we apply it only to the left of
$1$). 
For $z_0=1$, the repulsive periodic point is $z_1=0$, the
interval $[z_1,z_2]$ is $[0,\lambda^{-r}]$, and the growth factor is
$\frac{\alpha_0 \delta}{\delta}=2r$. 
For $z_0=1/2$, the repulsive
periodic point is $z_1=4$, the interval $[z_1,z_2]$ is
$[4-\frac{3}{2}\lambda^{-r},4]$, and the growth factor is
$\frac{2\delta}{\frac{3}{2}\lambda^{-r}}=3$.

Since $f_r^2(\lambda^{-r})=0$ and
$f_r(4-\frac{3}{2}\lambda^{-r})=y_1$, for  any non degenerate interval
$J\subset [0,1]\cup[y_1,4]$ there exists $k$ such that either
$|f_r^k(J)|\geq \frac{3}{2}|J|$ or  $f_r^k(J)$ contains one of the
points $0,4, y_1$.

\begin{lem}\label{lem:lem2.1}
If $J$ is a non degenerate subinterval of $[0,1]\cup[y_1,4]$, there
exists  $k\geq 0$ such that either $|f_r^k(J)|\geq \frac{3}{2}|J|$ or
$0\in f_r^k(J)$ or $4\in f_r^k(J)$ or $y_1\in \Int{f_r^k(J)}$.
\end{lem}

Since $f_r^2([0,\lambda^{-r}])=[0,4]$ and
$f_r^3([4-\frac{3}{2}\lambda^{-r},4])=f_r^2([y_1,4])= [0,4]$, applying
Lemma \ref{lem:point-repulsif} we obtain the next lemma.

\begin{lem}\label{lem:lem2.2}
If $J$ is a non degenerate subinterval containing either $0$ or $4$,
then there exists $k\geq 0$ such that $f_r^k(J)=[0,4]$.
\end{lem}

\demisaut\label{subsec:summary} The construction of $f_r\colon[0,4]\to
[0,4]$ is now finished. 
The map is $C^r$ on $[0,4]$ (and is
$C^{\infty}$ on $[0,4]\backslash\{ 1 \}$), and
$\|f_r'\|_{\infty}=\lambda^r$. 
Furthermore, if we put together Lemmas
\ref{lem:lem1.1}, \ref{lem:lem1.2}, \ref{lem:lem2.1} and
\ref{lem:lem2.2}, we see that for any non degenerate subinterval
$J\subset [0,4]$, there exists $k\geq 0$ such $f_r^k(J)=[0,4]$.

\begin{prop} \label{prop:summary}
$f_r \colon I\to I$ is $C^r$, mixing and $\|f_r'\|_{\infty}=\lambda^r$.
\end{prop}

\begin{rem}
If we identify the two end points $0$ and $4$, the map $f_r$ can be
seen as a mixing $C^r$ map on the torus, since
$f_r^{(k)}(0)=f_r^{(k)}(4)$ for every $k\geq 1$.
\end{rem}

\section{Markov chain associated with $f_r$} \label{sec:Markov}
We show that $f_r$ is a Markov map for a suitable countable partition.
The associated Markov chain reflects almost all topological properties
of the system $(I,f_r)$.

\subsection{Definition of the graph} \label{subsec:Markov-graph}
We explicit the Markov partition $V_r$ and the associated graph $G_r$.

Let $t^n_0=x_n<t^n_1<\cdots<t^n_{M_n}=y_n$ the local extrema of $f_r$
on $[x_n,y_n]$. 
Let
\begin{eqnarray*}
V_r &=& \{[t^n_{i-1},t^n_i]\colon 1\leq n, 1\leq i\leq M_n\} \\ &&
\cup\{[\lambda^{-kr}x_n,\lambda^{-kr}y_n]\colon 1\leq k\leq n\}\\ &&
\cup \{[\lambda^{-kr} y_{n+1},\lambda^{-kr}x_n]\colon 1\leq n,0\leq
k\leq n\}\\ && \cup\{[\lambda^{-nr}y_n,\lambda^{-(n-1)r}] \colon 2\leq
n\}\\ && \cup\{[\lambda^{-r}y_1,1/2],[1/2,1],[y_1,4]\}.
\end{eqnarray*}
The elements of $V_r$ have pairwise disjoint interior and their union
is $]0,4]$. 
We check that the map $f_r$ is monotone on each element of
$V_r$  and if $J\in V_r$ then $f_r(J)$ is a union of elements of
$V_r\cup\{0\}$.

\begin{itemize}
\item 
By Lemma \ref{lem:xn-yn}, $f_r$ is monotone on $[t_{i-1}^n,t_i^n]$,
$f_r([t_0^n,t_1^n])=[\lambda^{-nr}x_n,\lambda^{-nr}y_n]$ and
$f_r([t_{i-1}^n,t_i^n])=[\lambda^{-nr}y_{n+1},\lambda^{-nr}x_n]\cup
[\lambda^{-nr}x_n,\lambda^{-nr}y_n]$ if $2\leq i\leq M_n$.
\item 
By Lemmas \ref{lem:xn-yn} and \ref{lem:yn+1-xn}, $f_r$ is increasing
on $[y_{n+1},x_n]$ for all $n\geq 1$ and\\
$\disp
f_r([y_{n+1},x_n])=[\lambda^{-(n+1)r}y_{n+1},\lambda^{-nr}x_n]\\ 
= [\lambda^{-(n+1)r}y_{n+1},\lambda^{-nr}]\cup\bigcup_{k\geq n}
[\lambda^{-nr}x_{k+1},\lambda^{-nr}y_{k+1}]\cup
[\lambda^{-nr}y_{k+1},\lambda^{-nr}x_k].
$
\item
Since $f_r(x)=\lambda^r x$ for $x\in[0,\lambda^{-r}y_1]$ we have\\ 
-- $f_r([\lambda^{-kr}x_n,\lambda^{-kr}y_n])=[\lambda^{-(k-1)r} x_n,
\lambda^{-(k-1)r} y_n]$ for $1\leq k\leq n$ and this interval is an
element of $V_r$ except $\disp
[x_n,y_n]=\bigcup_{i=1}^{M_n}[t_{i-1}^n,t_i^n]$ which is a union of
elements of $V_r$.\\ 
-- $f_r([\lambda^{-kr}y_{n+1},\lambda^{-kr}x_n])=[\lambda^{-(k-1)r}y_{n+1},
\lambda^{-(k-1)r}x_n]$ for $1\leq k\leq n$.\\ 
--  $\disp
f_r([\lambda^{-(n+1)r}y_{n+1},\lambda^{-nr}])
=[\lambda^{-nr}y_{n+1},\lambda^{-(n-1)r}]$ \\
\hspace*{1em}$\disp =[\lambda^{-nr}y_{n+1},\lambda^{-nr}x_{n}] \cup
[\lambda^{-nr}x_{n},\lambda^{-nr}y_{n}] \cup
[\lambda^{-nr}y_{n},\lambda^{-(n-1)r}]$ for $n\geq 1$.
\item
$f_r$ is monotone on $[0,1/2]$, $[1/2,1]$ and $[y_1,4]$ (see
Subsection \ref{subsec:1.4}) and\\ 
-- $f_r([\lambda^{-r}y_1,1/2])=[y_1,4]$.\\ 
-- $\disp f_r([1/2,1])=[0,4]=\{0\}\cup\bigcup_{J\in V_r} J$.\\ 
-- $\disp f_r([y_1,4])=[\lambda^{-r}y_1,4]$\\
\hspace*{1em}$\disp =[\lambda^{-r}y_1,1/2]\cup[1/2,1]\cup [y_1,4]\cup
\bigcup_{n\geq 1} [y_{n+1},x_n]\cup \bigcup_{\mbox{\scriptsize
\begin{tabular}{c} $1\leq n$\\ $1\leq i\leq
M_n$\end{tabular}}}[t_{i-1}^n,t_i^n]$.
\end{itemize}

We define the directed graph $G_r$ as follows: the set of vertices of
$G_r$  is $V_r$  and there is an arrow from $J$ to $K$ if and only if
$K\subset f_r(J)$.  
The decomposition above of $f_r(J)$ into elements
of $V_r$ for all $J\in V_r$ gives an exhaustive description of the
arrows in $G_r$.

Notice that the graphs $G_r$ are identical for all $r\geq 1$.  
The only difference is the name of the vertices, corresponding to the
partition of $f_r$.

\subsection{Isomorphism between $f_r$ and the Markov chain}
Let $\Gamma_r^+$ be the set of all one-sided infinite sequences
$(D_n)_{n\geq 0}$ such that $D_n\in V_r$ and \makebox{$D_n\to
D_{n+1}$} $\forall n\in\IN$, and let $\Gamma_r$ be the set of all
two-sided infinite sequences $(D_n)_{n\in \IZ}$.  
We write $\sigma$
for the shift transformation in both spaces.  
$(\Gamma_r,\sigma)$ is
called the {\bf Markov chain associated with $f_r$}. 
As  the systems
$(\Gamma_r,\sigma)$ are isomorphic for all $r\geq 1$,  we just write
$(\Gamma,\sigma)$ when we want to talk about one of them without
referring to the partition associated with $f_r$.

We are going to build an isomorphism modulo countable sets between
$(I,f_r)$ and $(\Gamma_r^+,\sigma)$, that is a map $\phi_r \colon
I\backslash{\cal N}_r\longrightarrow \Gamma_r^+\backslash{\cal M}_r$
where ${\cal N}_r,\: {\cal M}_r$ are countable sets, $\phi_r$ is
bijective bimeasurable (in fact bicontinuous) and $\phi_r\circ
f_r=\sigma\circ\phi_r$.

\saut Define
\begin{eqnarray*}
{\cal P}_r&=& \{\lambda^{-kr}x_n,\lambda^{-kr}y_n\colon 1\leq k\leq
n\} \cup\{t^n_i\colon 1\leq n, 0\leq i\leq M_n\}\\ &
&\cup\{\lambda^{-nr}\colon 1\leq n\} \cup\{0,1/2,1,4\}
\end{eqnarray*}
and let ${\cal N}_r=\bigcup_{n\geq 0}f_r^{-n}({\cal P}_r)$  which is
countable.  
We have $f_r({\cal N}_r)={\cal N}_r$ and
$f_r(I\backslash{\cal N}_r)= I\backslash{\cal N}_r$.  
If $x\in
I\backslash{\cal P}_r$  then there is a unique $D\in V_r$ such that
$x\in D$ (in fact $x\in \Int{D}$). 
Hence if $x\in I\backslash{\cal
N}_r$, for  every $n\geq 0$ there is a unique $D_n\in V$ such that
$f_r^n(x)\in D_n$. 
Moreover $(D_n)_{n\geq 0}\in \Gamma_r^+$. 
We define
$$
\begin{array}{crcl} \phi_r \colon   & I\backslash{\cal N}_r
&\longrightarrow& \Gamma_r^+\\ 
&x&\mapsto& (D_n)_{n\geq 0}
\end{array}
$$
This application satisfies $\phi_r\circ f_r(x)= \sigma\circ\phi_r(x)$.

For any $(D_n)_{n\geq 0}\in \Gamma_r^+$, the set $J=\bigcap_{n\geq
0}f_r^{-n} (D_n)$ is a compact interval because $f_r$ is monotone on
each $D_n$.  
The map $f_r$ is mixing (Proposition \ref{prop:summary})
and $f_r^n(J)\subset D_n$,  hence $J$ is necessarily reduced to a
single point $\{x\}$. 
We define
$$
\begin{array}{crcl}
\psi_r \colon   &\Gamma_r^+ &\longrightarrow& I\\ &(D_n)_{n\geq 0}
 &\mapsto& x
\end{array}
$$
Let ${\cal M}_r=\psi_r^{-1}({\cal N}_r)$. 
The application $\psi_r$,
restricted to $\Gamma_r^+\backslash{\cal M}_r$, is the inverse of
$\phi_r$. 
Moreover, both  $\phi_r$ and $\psi_r$ are continuous.
Indeed, choose $x_0\in I\backslash{\cal N}_r$ and write $(D_n)_{n\geq
0}=\phi_r(x_0)$ and $J_n=\bigcap_{k=0}^n f_r^{-k}(D_k)$. 
The diameters
of the compact intervals $J_n$ tend to $0$, the point $x_0$ belongs to
$\Int{J_n}$ for every $n$, and for every  $x\in J_n\backslash {\cal
N}_r$ the sequence $\phi_r(x)$ begins with $(D_0,\cdots,D_n)$. 
Hence
$\phi_r$ is continuous. 
Inversely, fix $\disp \gamma_0=(D_n)_{n\geq
0}\in  \Gamma_r^+\backslash{\cal M}_r$, then for every sequence
$\gamma\in  \Gamma_r^+\backslash{\cal M}_r$  beginning with
$(D_0,\cdots,D_n)$ the point $\psi_r(\gamma)$ belongs to $J_n$ which
is an arbitrarily small neighbourhood  of $\psi_r(\gamma_0)$. 
Hence
$\psi_r$ is continuous too.

\saut Now, we are going to show that ${\cal M}_r$ is countable.  
It is
sufficient to show that $\psi_r^{-1}(x)$ is finite for any $x\in
{\cal N}_r$. 
For any $y\in I$ there are at most two elements of $V_r$
containing $y$.  
Let $x\in {\cal N}_r$. 
If there is a $k$ such that
$f_r^k(x)=0$ then  $\psi_r^{-1}(x)=\emptyset$. 
If there is a $k$ such
that $f_r^k(x)=4$ then $\psi^{-1}(x)$ is finite because $\psi^{-1}(4)$
contains only the constant sequence of symbol $[y_1,4]$.  
Otherwise
there exist $k,\,n$ such that $f_r^k(x)=x_n$ or $f_r^k(x)=y_n$. 
Thus
it is sufficient to focus on the points $x_n$ and  $y_n$.

We begin with $x_n$. 
The intervals $C_0=[y_{n+1},x_n]$ and
$D_0=[x_n,t^n_1]$ are the only two elements of $V_r$ containing
$x_n$. 
If we try to build $\disp (C_k)_{k\geq 0}$ and $\disp
(D_k)_{k\geq 0}$ which are elements of  $\psi_r^{-1}(x_n)$, we see
that there are only two possibilities, which are cycles, namely:
\begin{itemize}
\item
$C_0=[y_{n+1},x_n]\to
C_1=[\lambda^{-nr}y_{n+1},\lambda^{-nr}x_n]\to\cdots\to
C_{n+1}=C_0\to\cdots$
\item
$D_0=[x_n,t^n_1]\to
D_1=[\lambda^{-nr}x_n,\lambda^{-nr}y_n]\to\cdots\to
D_{n+1}=D_0\to\cdots$
\end{itemize}
Hence, $\Card{\psi_r^{-1}(x_n)}=2$.

The situation is the same for $y_n, n\geq 2$, with two slightly
different cycles, namely:
\begin{itemize}
\item
$C_0=[t^n_{M_n-1},y_n]\to
[\lambda^{-nr}x_n,\lambda^{-nr}y_n]\to\cdots\to C_{n+1}=C_0\to\cdots$
\item
$D_0=[y_n,x_{n-1}]\to [\lambda^{-nr}y_n,\lambda^{-(n-1)r}]\to
\cdots\to D_{n+1}=D_0\to\cdots$
\end{itemize}

A quick look at the map $f_r$ gives the last two cycles for $y_1$.

Consequently, $\Card{\psi_r^{-1}(x)}<+\infty$ for every $x\in {\cal
N}_r$, ${\cal M}_r$ is countable, and the map \makebox{$\phi_r \colon
I\backslash{\cal N}_r\longrightarrow \Gamma_r^+\backslash{\cal M}_r$}
is an isomorphism modulo countable sets.

\saut $\phi_r$ transforms any invariant measure that does not charge
${\cal N}_r$ into an invariant measure that does not charge ${\cal
M}_r$, and inversely.  
A measure supported by ${\cal N}_r$ or ${\cal
M}_r$ is of zero  entropy and the metric entropy $\mu\mapsto h_{\mu}$
is affine (see e.g. \cite{DGS}), thus
$h_{top}(f_r)=h(\Gamma_r^+,\sigma)$, where
$$
h(\Gamma_r^+,\sigma)=sup\{h_{\mu}\colon \mu\;
\sigma\mbox{-invariant } \mbox{measure on } \Gamma_r^+\},
$$ 
and
$\phi_r$ establishes a bijection between the sets of maximal measures.

On the other hand,  $h(\Gamma_r^+,\sigma)=h(\Gamma_r,\sigma)$ and
there is a bijection between the maximal measures of
$(\Gamma_r^+,\sigma)$ and those of $(\Gamma_r,\sigma)$, because the
latter is the natural extension of the former (see e.g.
\cite{Pet}). 
Recall that all systems  $(\Gamma_r,\sigma)$ are
identical and $(\Gamma,\sigma)$ represents equally one of them.  
Hence
the question of existence of maximal measure for $(I,f_r)$ can be
studied by looking at $(\Gamma,\sigma)$.

\begin{prop} \label{prop:isomorphism}
$h_{top}(f_r)=h(\Gamma,\sigma)$ and $(I,f_r)$ admits a maximal measure
if and only if $(\Gamma,\sigma)$ admits one.
\end{prop}

\subsection{Non existence of maximal measure} \label{subsec:non-existence}
Following the terminology of Vere-Jones \cite{Ver1} a transitive
Markov chain is either transient, positive recurrent or null
recurrent.  
According to a result of Gurevi\v{c} \cite{Gur2}, a
transitive Markov chain has a maximal measure if and only if  its
graph is positive recurrent. 
We do not give the definitions of transience, positive recurrence and 
null recurrence because we will only need a criterion due to Salama 
(Theorem 2.1(i) in \cite{Sal2}), which is stated below.

If $H$ is a (strongly) connected directed graph and
$(\Gamma_H,\sigma)$  is the associated Markov chain, i.e. the set of
all sequences $\disp (h_n)_{n\in\IZ}$ with $h_n\to h_{n+1}$ in $H$, we
define  $h(H)=h(\Gamma_H,\sigma)= \sup\{h_{\mu}\colon \mu\;
\sigma\mbox{-invariant probability on } \Gamma_H\}$.

\begin{theo} \label{theo:Gurevic}
{\bf (Gurevi\v{c})} Let $H$ be a connected directed graph and
$(\Gamma_H,\sigma)$  be the associated Markov chain. 
If its entropy
$h(H)$ is finite then  $(\Gamma_H,\sigma)$ admits a maximal measure if
and only if  $H$ is positive recurrent. 
In this case, the measure is
unique.
\end{theo}

\begin{theo} \label{theo:Salama}
{\bf (Salama)} Let $H$ be a connected directed graph.  
If there exists
a graph $H'$ such that $H\subsetnoteq H'$ and $h(H)=h(H')$ then $H$ is
transient.
\end{theo}

Next, we  compute $h(G_r)$ then we show that $G_r$ is transient, which
is enough to conclude that $f_r$ has no maximal measure by Proposition
\ref{prop:isomorphism}. 
As all graphs $G_r$ are identical,  it is
sufficient to focus on $G_1$.

\begin{prop}\label{prop:hG}
$h_{top}(f_r)=h(G_r)=\log\lambda$.
\end{prop}

\begin{demo}
It is already known that $h_{top}(f_r)=h(G_r)=h(G_1)$ by Proposition
\ref{prop:isomorphism}.

A subset $E\subset I$ is called $(n,\eps)$-separated for $f_1$ if for
any two distinct points $x,y$ in $E$ there exists $k$, $0\leq k <n$,
with $|f_1^k(x)-f_1^k(y)|>\eps$. 
Let $s_n(f_1,\eps)$ be the maximal
cardinality of an $(n,\eps)$-separated set. 
Then the topological
entropy of $f_1$ is given by the following formula (see
e.g. \cite{DGS}):
$$
h_{top}(f_1)=\lim_{\eps\to 0}\limsup_{n\to +\infty} \frac{1}{n}\log
s_n(f_1,\eps).
$$ 
Let $E$ be an $(n,\eps)$-separated set of $I$ of
maximal cardinality. 
As  $\|f_1'\|_{\infty}=\lambda$ (Proposition
\ref{prop:summary}),  we have $|f_1(x)-f_1(y)|\leq $
\makebox{$\lambda|x-y|$} for all $x,y\in I$. 
If $x,y$ are two distinct
points of $E$, there exists $k<n$ such that
$|f_1^k(x)-f_1^k(y)|>\eps$. 
But $|f_1^k(x)-f_1^k(y)|\leq
\lambda^n|x-y|$, hence $|x-y|\geq \lambda^{-n}\eps$ and
$$
\Card{E}=s_n(f_1,\eps)\leq \frac{\lambda^n}{\eps}+1.
$$
Consequently, $h_{top}(f_1)=h(G_1)\leq \log \lambda$.

\saut Now, let $H_n\subset G_1$ be the subgraph whose vertices are:
$$
\{[t^n_{i-1},t^n_i]\colon 1\leq i\leq M_n\}\cup\{[\lambda^{-k}x_n,
\lambda^{-k}y_n]\colon 1\leq k\leq n\}.
$$

The edges of $H_n$ are all possible edges of $G_1$  between two
vertices, namely:
\begin{itemize}
\item $[t^n_{i-1},t^n_i]\to [\lambda^{-n}x_n,\lambda^{-n}y_n]$ for
$1\leq i\leq M_n$,
\item $[\lambda^{-k}x_n,\lambda^{-k}y_n]\to [\lambda^{-k+1}x_n,
\lambda^{-k+1}y_n]$ for $2\leq k\leq n$,
\item $[\lambda^{-1}x_n,\lambda^{-1}y_n]\to [t^n_{i-1},t^n_i]$ for
$1\leq i\leq M_n$.
\end{itemize}
The graph $H_n$ is represented in Figure \ref{fig:Hn}.

\begin{figure}[htbp]
\begin{center}
\includegraphics{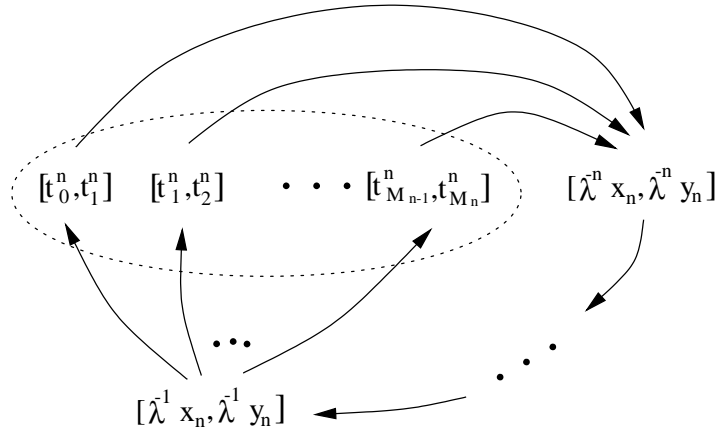}
\caption{the graph $H_n$; $\sigma^{n+1}$  is a full shift on the set
of vertices inside the dots.\label{fig:Hn}}
\end{center}
\end{figure}

The system $(H_n,\sigma^{n+1})$ is a full shift on $M_n$ symbols, plus
$n$ fixed points, thus $h(H_n,\sigma^{n+1})=\log M_n$ (see
e.g. \cite[p111]{DGS}) and $\disp h(H_n)=\frac{\log M_n}{n+1}$.

By definition of $M_n$,
$$
\lim_{n\to +\infty}\frac{\log M_n}{n+1}=\log \lambda.
$$
As $H_n$ is a subgraph of $G_1$,  $h(H_n)\leq h(G_1)$. 
Therefore $h(G_1)=\log\lambda$.
\end{demo}

\begin{prop} \label{prop:transient}
The graph $G_1$ is transient.
\end{prop}

\begin{demo}
We are going to build a Markov map $g$, very similar to  $f_1$, such
that $\|g'\|_{\infty}\leq \lambda$ and the Markov graph $H$
associated with $g$  expands strictly $G_1$.  
Suppose $g$ is already
built.  
The same argument as in the proof of Proposition \ref{prop:hG}
shows that $h(H)\leq \log\|g'\|_{\infty}\leq \log\lambda$. 
As $G_1\subset H$ we have $h(H)\geq h(G_1)$, thus $h(H)=h(G_1)$ by
Proposition \ref{prop:hG}.  
This is enough to conclude that $G_1$ is
transient by Theorem  \ref{theo:Salama}.

The map $g\colon I\to I$ is defined as $g(x)=f_1(x)$ for all $x\in
I\setminus [x_2,y_2]$. 
Let
$$
\widetilde{M}_2=M_2+2 \mbox{ and }
\widetilde{k}_2=\frac{2\lambda}{\widetilde{M}_2}
$$
and choose $\widetilde{s}_2\colon [0,\widetilde{M}_2]\to [-m_2,1]$
satisfying Properties (3)-(8) except that $M_2$ and $k_2$ are replaced
respectively by $\widetilde{M}_2$ and $\widetilde{k}_2$. 
Then we define $g$ on $[x_2,y_2]$ by
$$
g(x)=\lambda^{-2}\left[x_2+(y_2-x_2)\widetilde{s}_2\left(\widetilde{M}_2
\frac{x-x_2}{y_2-x_2}\right)\right].
$$  
By Properties (5) and (6),
$g'(x_2)=g'(y_2)=\lambda^{-2}\widetilde{M}_2\widetilde{k}_2=
2\lambda^{-1}$, thus $g'(x_2)=f_1'(x_2)$, $g'(y_2)=f_1'(y_2)$ and $g$
is $C^1$. 
Moreover for all $x\in [x_2,y_2]$,
$$
|g'(x)|\leq \lambda^{-2}
\widetilde{M}_2\|\widetilde{s}_2'\|_{\infty}\leq\lambda^{-1}\widetilde{M}_2
$$
thus $|g'(x)|<\lambda$  because
$\widetilde{M}_2=M_2+2=2\left[\frac{\lambda}{8}\right]+1<\lambda^2$.
Since $\|f_1'\|_{\infty}= \lambda$ by Proposition \ref{prop:summary},
one concludes that $\|g'\|_{\infty}\leq \lambda$

Define the Markov graph $H$ associated with $g$ as in Subsection
\ref{subsec:Markov-graph}, and denote by $W$ the set of vertices of $H$.
Compared to $V_1$, $W$ has two additional vertices
because $f_1$ has $M_2$ monotone pieces between $x_2$ and $y_2$ and
$g$ has $M_2+2$.  
If
$$
\widetilde{t}_i=x_2+\frac{i(y_2-x_2)}{M_2+2}
$$ 
for $0\leq i\leq M_2+2$ then it is not hard to check that the graph
$G_1$ is equal to $H$ deprived of the vertices
$[\widetilde{t}_{M_2},\widetilde{t}_{M_2+1}]$ and
$[\widetilde{t}_{M_2+1},\widetilde{t}_{M_2+2}]$ and all the edges that
begin or end at one of them. 
Consequently $G_1\subsetnoteq H$, which
ends the proof.
\end{demo}

\begin{rem}\label{rem:mu-n-dirac}
We can see intuitively what happens for an $f_r$-invariant measure
when its entropy tends to $\log\lambda$. 
On each finite subgraph
$H_n$, there is a measure of entropy $\frac{\log M_n}{n+1}$. 
This measure has a corresponding measure $\mu_n$ on the interval, the
support of which is contained in $\bigcup_{k=0}^n
[\lambda^{-kr}x_n,\lambda^{-kr}y_n]$ (in fact, the support of $\mu_n$
is exactly the Cantor set of all points which never escape from that
set).  
We have of course $h_{\mu_n}(f_r)\to \log\lambda$. 
But if we
consider what happens near $0$, we see that $\mu_n$ converges to
$\delta_0$, the Dirac measure at $0$, whose entropy is null.
\end{rem}

\section{Local entropy} \label{sec:calcul-entropie-locale}
We recall first some definitions due to Bowen \cite{Bow} and then  we
define the local entropy. 
There exist different definitions of local
entropy, we give here that of Buzzi \cite{Buz}.

\begin{defi}
Let $T\colon X\to X$ be a continuous map on a compact metric space
$X$.

The {\bf Bowen ball} of radius $r$ and  order $n$, centered at $x$ is
$B_n(x,r)=\{y\in X\colon d(T^k(y),T^k(x))<r,$ \makebox{$\forall\,
k=0,\cdots,n-1\}$}.

An {\bf $(\eps,n)$-separated set} of $Y\subset X$ is a subset
$E\subset Y$ such that $\forall\, y\not=y'$ in $E$,
\makebox{$\exists\; 0\leq k<n$},   $d(T^k(y),T^k(y'))>\eps$.  
The maximal cardinality of an $(\eps,n)$-separated set of $Y$ is denoted
by $s_n(T,\eps,Y)$.
\end{defi}

\begin{defi}
The {\bf local entropy} of $T$, $ h_{loc}(T)$, is defined as
$$
h_{loc}(T)=\lim_{\eps\to 0}\lim_{\delta\to
0}\limsup_{n\to\infty}\frac{1}{n} \sup_{x\in X}\log
s_n(T,\delta,B_n(x,\eps)).
$$
\end{defi}

\begin{rem}
An $(\eps,n)$-cover of $Y\subset X$ is a subset $S\subset X$ such that
$\disp Y\subset \bigcup_{x\in S} B_n(x,\eps)$.  
Some people use
$(\eps,n)$-covers instead of $(\eps,n)$-separated sets: it leads to
the same definition of the local entropy.
\end{rem}

Local entropy is interesting because it bounds the defect of upper
semicontinuity of the metric entropy $\mu\mapsto h_{\mu}(f)$.  
On a
compact Riemannian m-dimensional manifold, local entropy itself is
bounded by $\frac{m\log R(f)}{r}$,  where $R(f)$ is the spectral
radius of the differential and $r$ is the differential order.  
These
results are stated by Buzzi \cite{Buz} and follow works of Yomdin
\cite{Yom} and Newhouse \cite{New}. 
In particular, they directly imply
that a $C^{\infty}$ map on a compact Riemannian manifold always has a
maximal measure (this result can be found in Newhouse's  work
\cite{New}).  
These results are given in the next two theorems, the
second one is stated for interval maps only.

\begin{theo}
Let $T\colon X\to X$ be a continuous map on a compact metric space.
Assume that $\mu_n$ is a sequence of $T$-invariant measures on $X$,
converging to a measure $\mu$. 
Then
$$
 \limsup_{n\to\infty}h_{\mu_n}(T)\leq h_{\mu}(T)+h_{loc}(T).
$$
\end{theo}

\begin{theo} \label{theo:hloc}
Let $f\colon I\to I$ be a $C^r$ map on a compact interval $I$, $r\geq
1$, and let $\disp R(f)=\inf_{k\geq
1}\sqrt[k]{\|(f^k)'\|_{\infty}}$. 
Then the local entropy satisfies
$$
h_{loc}(f)\leq\frac{\log R(f)}{r}.
$$
\end{theo}

In our family of examples, the local entropy can be computed
explicitly.

\begin{prop} For every $n\geq 1$ 
the local entropy of $f_r$ is
$$
h_{loc}(f_r)=\frac{\log R(f_r)}{r}=\log \lambda.
$$
\end{prop}

\begin{demo}
The map $f_r$ is such that $\|f_r'\|_{\infty}\leq \lambda^r$
(Proposition \ref{prop:summary}) and $0$ is a fixed point with
$f_r'(0)=\lambda^r$. 
Hence $R(f_r)=\lambda^r$ and
$$
h_{loc}(f_r)\leq \frac{\log R(f_r)}{r}=\log\lambda
$$ 
according to Theorem \ref{theo:hloc}.

We are going to show the reverse inequality.

Fix $\eps>0$ and choose $n$ such that $\frac{1}{2n^2}<\eps$. 
Put $\delta_0=\frac{1}{2n^2M_n}$.  
If $x\in [x_n,y_n]$ satisfies
$f^{n+1}(x)\in [x_n,y_n]$ then $|f^i(x)-f^i(x_n)|<\eps$ for $0\leq
i\leq n+1$. 
We write $I_i= [t^n_{i-1},t^n_i]$ for $1\leq i\leq
M_n$. 
The length of each $I_i$ is $\delta_0$.

Choose a finite sequence $\omega=(\omega_0,\cdots,\omega_{p-1})$ with
$1\leq \omega_i\leq M_n$. 
Thanks to the isomorphism between $(I,f_r)$
and its Markov extension (Section \ref{sec:Markov}), there is a point
$x_{\omega}\in [x_n,y_n]$ with \makebox{$f^{(n+1)i}(x_{\omega})\in
I_{\omega_i}$} for $0\leq i\leq p-1$. 
Consider the set
$E_{n,p}=\{x_{\omega}\colon \omega=(\omega_0, \cdots,\omega_{p-1}),
\omega_i \mbox{ odd}\}$.  
The cardinality of $E_{n,p}$ is
$$
\left(\frac{M_n+1}{2}\right)^p\geq
\left(\frac{\lambda^n}{4n^2}\right)^p
$$
by Lemma \ref{lem:lambda-n} (ii).  
If $x \in E_{n,p}$ then
$|f^k(x_n)-f^k(x)|<\eps$ for $0\leq k<(n+1)p$.  
Moreover, if
$x_{\omega},x_{\omega'}$ are two distinct elements of $E_{n,p}$, then
there exists $0\leq i\leq p-1$ with $|\omega_i-\omega'_i| \geq 2$,
hence $|f^{(n+1)i}(x_{\omega})-f^{(n+1)i}(x_{\omega'})|\geq
\delta_0$. 
Consequently, $E_{n,p}$ is an
$((n+1)p,\delta,B_{(n+1)p}(x_n,\eps) )$-separated set for every
$\delta< \delta_0$, and
$$
h_{loc}(f_r)\geq \lim_{n\to +\infty}\limsup_{p\to\infty} \frac{\log
\Card{E_{n,p}}}{(n+1)p}\geq \log\lambda.
$$
\end{demo}

This computation shows that the bound $\frac{\log R(f)}{r}$ is a sharp
one to estimate the local entropy. 
Moreover, we remarked (Remark
\ref{rem:mu-n-dirac})  that there exists a sequence of measures
$\mu_n$ converging to the Dirac measure $\delta_0$, with
$h_{\mu_n}(f_r)\to h_{top}(f_r)$. 
Hence, the local entropy is exactly
the defect of upper semicontinuity of the metric entropy in this case.


\noindent
S. {\sc Ruette} - Institut de Math\'ematiques de Luminy - CNRS - case
907 - 163 avenue de Luminy - 13288 Marseille cedex 9 - France -  {\it
e-mail :} {\tt  ruette@iml.univ-mrs.fr}
 

\begin{thebibliography}{10}

\bibitem{ALM}
Ll. Alsed{\`a}, J.~Llibre, and M.~Misiurewicz.
\newblock {\em Combinatorial dynamics and entropy in dimension one}.
\newblock World Scientific Publishing Co. Inc., River Edge, NJ, 1993.

\bibitem{BCop}
L.~S. Block and W.~A. Coppel.
\newblock {\em Dynamics in One Dimension}.
\newblock Lecture Notes in Mathematics 1513. Springer-Verlag, 1992.

\bibitem{Blo3}
A.~M. Blokh.
\newblock On sensitive mappings of the interval.
\newblock {\em Russ. Math. Surv.}, 37:203--204, 1982.

\bibitem{Blo1}
A.~M. Blokh.
\newblock Decomposition of dynamical systems on an interval.
\newblock {\em Russ. Math. Surv.}, 38:133--134, 1983.

\bibitem{Bow}
R.~Bowen.
\newblock Entropy for group endomorphisms and homogeneous spaces.
\newblock {\em Trans. Amer. Math. Soc.}, 153:401--414, 1971.

\bibitem{Bow2}
R.~Bowen.
\newblock Periodic points and measures for {A}xiom ${A}$ diffeomorphisms.
\newblock {\em Trans. Amer. Math. Soc.}, 154:377--397, 1971.

\bibitem{Bow3}
R.~Bowen.
\newblock Some systems with unique equilibrium states.
\newblock {\em Math. Syst. Theory}, 8:193--202, 1974.

\bibitem{Buz4}
J.~Buzzi.
\newblock personal communication.

\bibitem{Buz}
J.~Buzzi.
\newblock Intrinsic ergodicity of smooth interval maps.
\newblock {\em Israel J. Math.}, 100:125--161, 1997.

\bibitem{Buz3}
J.~Buzzi.
\newblock Specification on the interval.
\newblock {\em Trans. Amer. Math. Soc.}, 349(7):2737--2754, 1997.

\bibitem{DGS}
M.~Denker, C.~Grillenberger, and K.~Sigmund.
\newblock {\em Ergodic Theory on Compact Spaces}.
\newblock Lecture Notes in Mathematics 527. Springer-Verlag, 1976.

\bibitem{GZ}
B.~M. Gurevich and A.~S. Zargaryan.
\newblock A continuous one-dimensional mapping without a measure with maximal
  entropy.
\newblock {\em Functional Anal. Appl.}, 20(no. 2):134--136, 1986.

\bibitem{Gur1}
B.~M. Gurevi\v{c}.
\newblock Topological entropy of enumerable {M}arkov chains.
\newblock {\em Soviet. Math. Dokl.}, 10(no. 4):911--915, 1969.

\bibitem{Gur2}
B.~M. Gurevi\v{c}.
\newblock Shift entropy and {M}arkov measures in the path space of a
  denumerable graph.
\newblock {\em Soviet. Math. Dokl.}, 11:744--747, 1970.

\bibitem{Hof}
F.~Hofbauer.
\newblock On intrinsic ergodicity of piecewise monotonic transformations with
  positive entropy.
\newblock {\em Israel J. Math.}, {\bf I} 34:213--237, 1979.
\newblock {\bf II} 38:107--115, 1981.

\bibitem{Hof2}
F.~Hofbauer.
\newblock Piecewise invertible dynamical systems.
\newblock {\em Probab. Th. Rel. Fields}, 72:359--386, 1986.

\bibitem{MS2}
M.~Misiurewicz and W.~Szlenk.
\newblock Entropy of piecewise monotone mappings.
\newblock {\em Studia Math.}, 67(1):45--63, 1980.

\bibitem{New}
S.~E. Newhouse.
\newblock Continuity properties of entropy.
\newblock {\em Annals of Mathematics}, 129:215--235, 1989.

\bibitem{Pet}
K.~Petersen.
\newblock {\em Ergodic theory}.
\newblock Cambridge University Press, 1983.

\bibitem{Rue2}
D.~Ruelle.
\newblock Thermodynamic formalism for maps satisfying positive expansiveness
  and specification.
\newblock {\em Nonlinearity}, 5(6):1223--1236, 1992.

\bibitem{Sal}
I.~A. Salama.
\newblock {\em Topological entropy and classification of countable chains}.
\newblock PhD thesis, University of North Carolina, Chapel Hill, 1984.

\bibitem{Sal2}
I.~A. Salama.
\newblock Topological entropy and recurrence of countable chains.
\newblock {\em Pacific J. Math.}, 134(no. 2):325--341, 1988.
\newblock Errata, {\it Pacific J. Math.}, 140(no. 2):397, 1989.

\bibitem{Ver1}
D.~Vere-Jones.
\newblock Geometric ergodicity in denumerable {M}arkov chains.
\newblock {\em Quarterly J. Math.}, 13:7--28, 1962.

\bibitem{Yom}
Y.~Yomdin.
\newblock Volume growth and entropy.
\newblock {\em Israel J. Math.}, 57:285--300, 1987.
\newblock $C^k$-resolution of semi-algebraic mappings - Addendum to ``Volume
  growth and entropy'', {\it Israel J. Math.}, 57:301-318, 1987.

\end{thebibliography}
\end{document}